\documentclass[10pt]{article}
\usepackage{amsmath, amsthm, amssymb, amsfonts, mathrsfs, mathtools, upgreek}
\usepackage{graphicx}
\usepackage{makeidx}
\usepackage{tikz-network}
\usepackage[left=2cm,right=2cm,top=2cm,bottom=2cm]{geometry}
\usepackage{caption}
\usepackage{subcaption}
\usepackage[section]{placeins}
\usepackage{enumitem}
\usepackage{tikz}
\usepackage{stmaryrd}
\usepackage{upquote}
\usepackage{authblk}

\usepackage{multicol}
\usepackage{booktabs}
\usepackage{hyperref}
\hypersetup{colorlinks=true, citecolor={blue}, linkcolor=blue, filecolor=magenta, urlcolor=cyan}
\usetikzlibrary{arrows.meta, positioning, decorations.pathreplacing}
\usepackage{cleveref}
\newtheorem{theorem}{Theorem}[section]

\newtheorem{lemma}[theorem]{Lemma}

\theoremstyle{definition}
\newtheorem*{defi}{Definition}

\numberwithin{equation}{section}
\newcommand{\divides}{\mid}

\newcommand{\tr}{\operatorname{tr}}

\newcommand{\cay}{\operatorname{Cay}}

\newcommand{\iu}{\mathbf{i}}

\newcommand{\spec}{\operatorname{Spec}}
\newcommand{\lcm}{\operatorname{lcm}}
\newcommand{\trace}{\operatorname{trace}}
\newcommand{\rank}{\operatorname{rank}}
\newcommand{\ord}{\operatorname{ord}}
\newcommand{\Nl}{\mathbb{N}}
\newcommand{\Zl}{\mathbb{Z}}

\newcommand{\Ql}{\mathbb{Q}}

\newcommand{\Cl}{\mathbb{C}}

\newcommand{\eu}{\mathbf{e}_u}
\newcommand{\ev}{\mathbf{e}_v}
\newcommand{\ex}{\mathbf{e}_x}

\newcommand{\ld}{\lambda}
\newcommand{\ee}{\mathbf{e}}
\newcommand{\mh}{\mathcal{M}_h}
\usepackage{titling}
\title{\textbf{A complete classification of perfect state transfer in Grover walks on abelian Cayley graphs}}

\author{ Koushik Bhakta and Bikash Bhattacharjya\\
	Department of Mathematics\\
	Indian Institute of Technology Guwahati, India\\
	b.koushik@iitg.ac.in, b.bikash@iitg.ac.in }

\date{}
\begin{document}
\maketitle
	
	\vspace{-1cm}
	
	\begin{center}{\textbf{Abstract}}\end{center}

	\noindent
	We completely classify connected Cayley graphs over finite abelian groups that exhibit perfect state transfer in Grover walks. First, we show that a Cayley graph over any finite group exhibiting perfect state transfer has valency $1$, $2$, or $4$. For valency $4$, we prove that an abelian Cayley graph $\operatorname{Cay}(\Gamma,S)$ exhibits perfect state transfer if and only if $S=\{a,-a,a+h,-a+h\}$, where $h$ has order $2$, and $4\nmid |\Gamma|/2$. Combining this characterization with the valency-one and valency-two cases, we explicitly determine all pairs  $(\Gamma,S)$ for which $\operatorname{Cay}(\Gamma,S)$ exhibits perfect state transfer. Our results extend the classification of circulant graphs of valency at most four obtained by Kubota and Yoshino [Circulant graphs with valency up to 4 that admit perfect state transfer in Grover walks, \emph{J. Comb. Theory Ser. A.}, 216 (2025) 106064] to Cayley graphs over arbitrary finite abelian groups and to all valencies.
	
	\vspace*{0.2cm}
	
	
	
	\section{Introduction}
	This article studies a problem in algebraic graph theory motivated by quantum computation and quantum information. We focus on the problem of transferring a quantum state from one vertex to another with probability $1$ in the Grover walk~\cite{hig3}, a widely studied model of discrete-time quantum walks. Quantum walks are useful tools in quantum computation and quantum information, with applications to quantum algorithms~\cite{quantum_algorithm}, quantum simulation~\cite{simulation}, and quantum cryptography~\cite{cryptography}. According to the time evolution, quantum walks are broadly divided into two types: continuous-time quantum walks~\cite{state} and discrete-time quantum walks~\cite{godsil_dqw}. The transfer of a quantum state with probability $1$ is known as \emph{perfect state transfer}~\cite{state}. Perfect state transfer has been extensively studied for continuous-time quantum walks (see \cite{sofia_pst,guo_ctw,peter,singh,bi_cayley,quasi_abelian,cao_abelian}), whereas comparatively less is known in the discrete-time setting. We refer to~\cite{guo_dqw,kubota2, zhan1} for perfect state transfer in discrete-time quantum walks. For other dynamical phenomena in discrete-time quantum walks, we refer to~\cite{banerjee2,bethetree,sarkar,yoshie2,zhan3}. In this paper, we completely classify connected Cayley graphs over finite abelian groups that exhibit perfect state transfer in Grover walks (see Theorem~\ref{complete_classification}). To the best of our knowledge, a complete classification of perfect state transfer on abelian Cayley graphs has not previously been studied in either the discrete-time or continuous-time quantum walk setting. We first recall the Grover walk and the notion of perfect state transfer, and then present motivation for our work and our main results.
	
	\subsection{Grover walks and perfect state transfer}
	
	Throughout the paper, all graphs are finite and simple, with no isolated vertices. Let $G:=(V(G),E(G))$ be a graph with vertex set $V(G)$ and edge set $E(G)$. We write the elements of $E(G)$ as $uv$, where $u,v\in V(G)$, $u\neq v$, and $uv$ and $vu$ represent the same edge. For each edge $uv\in E(G)$, the ordered pairs $(u,v)$ and $(v,u)$ are called the \emph{arcs} associated with $uv$. We define $\mathcal{A}(G):=\{(u, v), (v, u):uv\in E(G) \}$, the set of all symmetric arcs of $G$. For an arc $a:=(u,v)$, the vertices $u$ and $v$ are called the \emph{origin} and \emph{terminus} of $a$, respectively, and we write $o(a)=u$ and $t(a)=v$.  The inverse arc of $a$, denoted $a^{-1}$, is the arc  $(v,u)$.
	
	The \emph{boundary matrix} $N:=N(G)\in \mathbb{C}^{V(G)\times \mathcal{A}(G)}$ of $G$ is defined by $N_{u,a}=\left(1/\sqrt{\deg u}\right)\delta _{u, t(a)},$ where $\delta_{a,b}$ is the Kronecker delta function. The \emph{shift matrix} $R:=R(G)\in \mathbb{C}^{\mathcal{A}(G) \times \mathcal{A}(G)}$ of $G$ is defined by $R_{a,b}=\delta_{a,b^{-1}}$. Define the \emph{time evolution matrix} $U:=U(G)\in \mathbb{C}^{\mathcal{A}(G)\times \mathcal{A}(G)}$ of $G$ by $$U=R(2N^*N-I).$$ The time evolution matrix is also known as the Grover transition matrix. A discrete-time quantum walk on a graph G is determined by a unitary matrix, which acts on the complex functions of the symmetric arcs of G. The discrete-time quantum walks defined by $U$ are called the \emph{Grover walks}. The Grover walks are also known as a special case of bipartite walks~\cite{chen}.
	
	For $\Phi\in \Cl^{\mathcal{A}(G)}$, we denote by $\|\Phi\|$ the Euclidean norm of $\Phi$.  A vector $\Phi\in \Cl^{\mathcal{A}(G)}$ is said to be a \emph{state} if $\|\Phi\|=1$. We say that \emph{perfect state transfer} occurs from a state $\Phi$ to another state $\Psi$ at time $\tau\in \Nl$ if there exists a unimodular complex number $\gamma$ such that $U^\tau\Phi=\gamma \Psi$. 
	
	We focus on the transfer of states localized at individual vertices of a graph. A state $\Phi\in \Cl^{\mathcal{A}(G)}$ is called \emph{vertex-type} in $G$ if there exists a vertex $u\in V(G)$ such that $\Phi=N^* \eu$. We denote the vertex-type state corresponding to the vertex $u$ by $\Phi_u:=N^*\eu$. Throughout this paper, we restrict our attention to vertex-type states. For the motivation behind vertex-type states, we refer the reader to \cite{kubota1}.
	\begin{defi}
		A graph exhibits \emph{perfect state transfer} from vertex $u$ to another vertex $v$ at time $\tau$ if there exists a complex number $\gamma$ with $|\gamma|=1$ such that $U^\tau \Phi_u=\gamma \Phi_v$.
	\end{defi}
	Throughout this paper, perfect state transfer refers to transfer from a vertex $u$ to a distinct vertex $v$, that is $v\ne u$. We say that a graph $G$ exhibits {perfect state transfer} if there exist vertices $u$ and $v$ in $G$ and a positive integer $\tau$ such that perfect state transfer occurs from $u$ to $v$ at time $\tau$.
	\subsection{Motivation and main results}
	
	Let $\Gamma$ be a finite group with identity element $1$, and let $S\subseteq\Gamma\setminus\{1\}$ satisfy $S\neq\emptyset$ and $S=S^{-1}$. The \emph{Cayley graph} $\cay(\Gamma,S)$ has vertex set $\Gamma$, where two vertices $u$ and $v$ are adjacent if $uv^{-1}\in S$. When $\Gamma$ is abelian, we use additive notation for the group operation. Then the identity element is denoted by $0$. A particularly important class of abelian Cayley graphs is \emph{circulant graphs}, namely, Cayley graphs on cyclic groups.

	In 2019, Zhan~\cite{zhan1} discovered an infinite family of circulant graphs with valency $4$ exhibiting perfect state transfer. This provided one of the first substantial families of circulant graphs with perfect state transfer in Grover walks. Later, in 2025, Kubota and Yoshino~\cite{kubota2} extended this work and obtained a complete characterization of circulant graphs of valency at most $4$ that exhibit perfect state transfer. The cases of valency $1$ and $2$ are relatively simple: a connected circulant graph of valency $1$ is $K_2$, while a connected circulant graph of valency $2$ is a cycle. Thus, the main difficulty lies in valencies $3$ and $4$. Kubota and Yoshino proved that no circulant graph with valency $3$ exhibits perfect state transfer, whereas for valency $4$ they obtained a complete classification; see, for instance, Theorem~\ref{kubota_2_main}. Their classification recovers Zhan's infinite family, gives another infinite family, and rules out perfect state transfer in all other cases.
	\begin{theorem}[{\cite[Theorem~1.3]{kubota2}}]\label{kubota_2_main}
		Let $G:=\cay(\Zl_{2\ell},\{\pm a,\pm b\})$ be a connected $4$-regular circulant graph, where $1\le a<b\le\ell-1$. Then $G$ exhibits perfect state transfer from $u$ to $v$ if and only if $a+b=\ell$, $v=u+\ell$, and $4\nmid \ell$.
	\end{theorem}
	
	These results naturally raise the question of what happens beyond valency $4$. In particular, it remained open whether circulant graphs of larger valency can exhibit perfect state transfer. We resolve this question in a stronger form. We prove that no Cayley graph on a finite group, not necessarily abelian, with valency greater than $4$ exhibits perfect state transfer. We then turn to the remaining nontrivial case of valency $4$ and completely characterize perfect state transfer on Cayley graphs over finite abelian groups. Thus, our results extend the classification from circulant graphs to the broader class of abelian Cayley graphs and, at the same time, establish a general upper bound on the valency of any Cayley graph exhibiting perfect state transfer. Our first main result is the following.

	\begin{theorem}\label{main1}
		Let $\Gamma$ be a finite group, and let $\emptyset\neq S\subseteq\Gamma\setminus\{1\}$ satisfy $S=S^{-1}$. If  $\cay(\Gamma,S)$ exhibits perfect state transfer, then $|S|\in\{1,2,4\}$.
	\end{theorem}
	
	The classification of circulant graphs with valency $4$ obtained by Kubota and Yoshino requires considerable effort, particularly in ruling out cases that do not exhibit perfect state transfer. An important part of their argument relies on arithmetic properties of the rings of integers of cyclotomic fields. In the more general setting of abelian Cayley graphs, the main difficulty is again the necessity: assuming perfect state transfer, why must the connection set have a very specific form?
	
	Our approach is fundamentally different and introduces new techniques to the study of perfect state transfer in quantum walks. We use elementary linear algebra, the structure of finite abelian groups, and basic facts about algebraic integers. These tools allow us to show that perfect state transfer strongly restricts the form of the connection set. This gives a simple and unified characterization, and extends the result of Kubota and Yoshino from circulant graphs to Cayley graphs on arbitrary finite abelian groups. Our second main result is the following.

	\begin{theorem}\label{main2}
		Let $\Gamma$ be a finite abelian group, and let $S\subseteq\Gamma\setminus\{0\}$ satisfy
		$S=-S$, $\langle S\rangle=\Gamma$ and $|S|=4$. Then $\operatorname{Cay}(\Gamma,S)$ exhibits
		perfect state transfer if and only if there exist $a,h\in\Gamma$ with $h$ of order $2$ such that $4$ does not divide $|\Gamma|/2$ and  $S=\{a,\,-a,\,a+h,\,-a+h\}$. Moreover, if perfect state transfer occurs, then it occurs from each vertex $u$ only to $u+h$ and the minimum transfer time is $\lcm(2,\ |\Gamma|/2)$.
	\end{theorem}

	The connected abelian Cayley graphs of valency $1$ and $2$ are, respectively, $K_2$ and cycles, for which perfect state transfer is already completely understood. This together with Theorem~\ref{main2} gives a complete classification of the underlying graphs that exhibit perfect state transfer. However, classifying the graphs alone does not determine all possible pairs $(\Gamma,S)$. For example, The graphs $\cay(\Zl_4,\{1,-1\})$ and $\cay\bigl(\Zl_2\times\Zl_2,\{(1,0),(0,1)\}\bigr)$
	are both the cycle $C_4$. However, the groups are not isomorphic. Thus, to obtain a complete classification of the Cayley representations, we also need to determine the possible groups and connection sets. Two pairs $(\Gamma,S)$ and $(\Gamma',S')$ are called  \emph{equivalent} if there exists a group isomorphism $\varphi:\Gamma\longrightarrow\Gamma'$
	such that $\varphi(S)=S'$. Such an isomorphism also induces a graph isomorphism between the corresponding Cayley graphs.  Our third and final main theorem gives, up to such equivalence, all pairs  $(\Gamma,S)$ for which $\cay(\Gamma,S)$ exhibits perfect state transfer.
	
	\begin{theorem}\label{complete_classification}
		Let $\Gamma$ be a finite abelian group, and let $S\subseteq\Gamma\setminus\{0\}$ satisfy $S=-S$ and $\langle S\rangle=\Gamma$. Then $\cay(\Gamma,S)$ exhibits perfect state transfer if and only if $(\Gamma,S)$ is equivalent to one of the following pairs.
		\begin{enumerate}[label=\textup{(\roman*)}]
			\item $\bigl(\Zl_2,\{1\}\bigr)$;	
			\item $\bigl(\Zl_2\times\Zl_2,\{(1,0),(0,1)\}\bigr)$;
			\item $\bigl(\Zl_{2n},\{1,-1\}\bigr)$, where $n\geq2$;
			\item $\bigl(\Zl_{2m},\{1,-1,m+1,m-1\}\bigr)$, where
			$m\geq3$ and $4$ does not divide $m$; or
			\item $\bigl(\Zl_m\times\Zl_2,
			\{(1,0),(-1,0),(1,1),(-1,1)\}\bigr)$, where
			$m\geq6$ and $m\equiv2\pmod4$.
		\end{enumerate}
	\end{theorem}
	%	
	%	\begin{corollary}\label{main3}
		%		Let $G$ be a connected Cayley graph over a finite abelian group. Then $G$ exhibits perfect state transfer if and only if $G$ is isomorphic to $K_2$, cycle $C_{2n}$, $n\ge2$, or $C_m[\overline{K}_2]$, where $m\ge3$ and $4\nmid m$.
		%	\end{corollary}
	%==================================

	% The remainder of the paper is organized as follows. In Section~2, we give necessary preliminaries and establish spectral criteria for perfect state transfer in abelian Cayley graphs. In Section~3, we prove Theorem~\ref{main1}, which restricts the possible valencies of Cayley graphs exhibiting perfect state transfer. Section~4 is devoted to abelian Cayley graphs of valency $4$, where we prove Theorem~\ref{main2}. Finally, in Section~5, we combine these results to obtain the complete classification stated in Corollary~\ref{main3}.

	\section{Preliminary results}
	
	In this section, we introduce some matrices and characterizations of perfect state transfer that will be used in subsequent sections. Let $G=(V(G),E(G))$ be a finite simple graph. The \emph{adjacency matrix} $A:=A(G)\in\mathbb C^{V(G)\times V(G)}$ of $G$ is defined by
	\[
	A_{u,v}=
	\begin{cases}
		1 & \text{if }uv\in E(G),\\
		0 & \text{otherwise}.
	\end{cases}
	\]
	The \emph{discriminant matrix} of $G$ is defined by $P:=P(G)=NRN^*$. If  $G$ is a regular graph, then the matrix $P$ can be expressed in terms of its adjacency matrix $A$.
	\begin{lemma}[{\cite[Lemma~3.1]{bhakta4}}]\label{reg}
		If $G$ is a $k$-regular graph, then $P=\frac{1}{k}A$. Further, $1$ is an eigenvalue of $P$, and for any eigenvalue $\mu$ of $P$,  $|\mu|\leq1$.
	\end{lemma}
	Kubota and Segawa~\cite{kubota1} studied perfect state	transfer in Grover walks using Chebyshev polynomials of the first kind. The \emph{Chebyshev polynomials of the first kind}, denoted by $T_n(x)$, are defined by the recurrence relation $T_0(x)=1$, $T_1(x)=x$, and $T_n(x)=2xT_{n-1}(x)-T_{n-2}(x)$ for $n\geq2$. It is well known that $T_n(\cos\theta)=\cos (n\theta)$. Consequently, for $n\geq1$ and $x\in[-1,1]$, we have $|T_n(x)|=1$ if and only if $x=\cos(j\pi/n)$ for some integer $j$.
	\begin{lemma}[{\cite[Lemma~3.1]{kubota1}}]\label{P_and_U}
		Let $G$ be a finite simple graph with the transition matrix $U$ and the discriminant matrix $P$. Then $NU^kN^*=T_k(P)$ for any nonnegative integer $k$.
	\end{lemma}
	Lemma~\ref{P_and_U} can be proved by induction.	This identity provides a connection between Grover walks and classical random walks via the discriminant matrix $P$.
	\subsection{Characterizations of Perfect State Transfer}
	
	Kubota and Yoshino~\cite{kubota2} showed that a graph $G$ exhibits perfect state transfer from $u$ to $v$	at time $\tau$ if and only if $T_\tau(P)\mathbf e_u=\gamma\mathbf e_v$ for some $\gamma\in\{-1,1\}$. Later, Guo and Schmeits~\cite{guo_sch}
	showed that necessarily $\gamma=1$. For completeness, we provide a proof of the following characterization.
	\begin{lemma}\label{pst_vertex}
		Let $u$ and $v$ be two distinct vertices of a graph $G$. Then $G$ exhibits perfect state transfer from $u$ to $v$ at time $\tau$ if and only if $T_\tau(P)\eu=\ev$.
	\end{lemma}
	
	\begin{proof}
		Suppose $G$ exhibits perfect state transfer from $u$ to $v$ at time $\tau$. Then $U^\tau N^*\eu=\gamma N^*\ev$.
		Applying $N$ to both sides and using Lemma~\ref{P_and_U}, we obtain $T_\tau(P)\eu=\gamma\ev$. 
		
		Let $w=(\sqrt{\deg x})_{x\in V(G)}$. Then $Pw=w$, giving that $T_\tau(P)w=w$. This, along with $T_\tau(P)e_u=\gamma e_v$, give that  $\sqrt{\deg u}=w^tT_\tau(P)e_u=\gamma\sqrt{\deg v}$.
		Since $\deg u$ and $\deg v$ are positive, and $|\gamma|=1$, we find that $\gamma=1$ 
		
		Conversely, suppose that
		$T_\tau(P)\eu=\ev$. By Lemma~\ref{P_and_U}, it follows that
		\[\langle N^*\ev,U^\tau N^*\eu\rangle=\ev^*NU^\tau N^*\eu=\ev^*T_\tau(P)\eu=1.\] 
		Since $NN^*=I$ and $U$ is unitary, both $N^*\ev$ and $U^\tau N^*\eu$ are unit vectors. Consequently,
		\[
		\|U^\tau N^*\eu-N^*\ev\|^2
		=2-2\Re\langle N^*\ev,
		U^\tau N^*\eu\rangle=0.
		\]
		Therefore $	U^\tau N^*\eu=N^*\ev$. Hence perfect state transfer occurs from $u$ to $v$.
	\end{proof}
	
	For Cayley graphs, vertex transitivity gives a stronger characterization. Since $T_\tau(P)$ commutes with every graph automorphism, perfect state transfer between two vertices determines its action on all vertices. This leads to the following result.
	\begin{lemma}\label{pst_cayley}
		Let $G:=\cay(\Gamma,S)$ be a Cayley graph over a finite group $\Gamma$. Then $G$ exhibits perfect state transfer at time $\tau$ if and only if  $T_\tau(P)$ is a permutation matrix of order $2$ with no fixed points.
	\end{lemma}
	\begin{proof}
		Suppose $G$ exhibits perfect state transfer from the vertex  $u$ to another vertex $v$ at time $\tau$. By Lemma~\ref{pst_vertex},  we have $ T_{\tau}(P)\eu=\ev$.  Using Lemma~\ref{reg}, $T_\tau(P)$ is a polynomial in $A$. Further, $\mathbf{e}^t_{xg}A^k\mathbf{e}_{yg} =\mathbf{e}^t_{x}A^k\mathbf{e}_{y}$ for any $x,y,g\in\Gamma$ and nonnegative integer $k$. Thus it follows that $\mathbf{e}^t_{xg}T_\tau(P)\mathbf{e}_{yg} =\mathbf{e}^t_{x}T_\tau(P)\mathbf{e}_{y}$. Combining this with $T_{\tau}(P)\eu=\ev$ and the symmetricity of $T_\tau(P)$, we deduce that $T_\tau(P)$ is a permutation matrix of order $2$ with no fixed points. The converse part follows from Lemma~\ref{pst_vertex}.
	\end{proof}
	\subsection{Algebraic integers and trace inequalities}
	We recall a few basic notions that will be used	throughout the subsequent sections. A complex number $\alpha$ is an \emph{algebraic number} if it is a 	root of a nonzero polynomial in $\mathbb{Q}[x]$. The unique monic irreducible polynomial $m_\alpha(x)\in\Ql[x]$ having $\alpha$ as a root is called the \emph{minimal polynomial} of $\alpha$ over $\mathbb{Q}$.  The roots of the minimal polynomial $m_\alpha(x)$ are called the \emph{algebraic conjugates} of $\alpha$ over $\mathbb{Q}$. An algebraic number $\alpha$ is called an \emph{algebraic integer} if it is a root of a monic polynomial with coefficients in $\mathbb{Z}$. It is known that the set of all algebraic integers forms a commutative ring, and that a rational number is an algebraic integer if and only if it is an integer.
	\begin{lemma}\label{algebraic_pst}
		Let $\mu$ be an eigenvalue of the discriminant of a Cayley graph $G$. If $G$ exhibits perfect state transfer, then $2\mu$ is an algebraic integer.
	\end{lemma}
	\begin{proof}
		Suppose $G$ exhibits perfect state transfer at time $\tau$. Then by Lemma~\ref{pst_cayley}, $T_\tau(\mu)=\pm1$. Hence $\mu=\cos(j\pi/\tau)$ for some integer $j$. This implies that $2\mu=\omega_{2\tau}^j+\omega_{2\tau}^{-j}$, where $\omega_{n}:=\exp(2\pi \iu/n)$. Hence $2\mu$ is an algebraic integer.
	\end{proof}
	\begin{lemma}\label{alg_kro}
		If $\alpha$ is a real algebraic integer with every algebraic conjugate of $\alpha$ has absolute value at most $1$, then $\alpha\in\{-1,0,1\}$.
	\end{lemma}
	\begin{proof}
		Suppose $\alpha\neq 0$, and let $p(x)=x^n+a_{n-1}x^{n-1}+\cdots+a_0\in\Zl[x]$   be the minimal polynomial of $\alpha$. Let
		$\alpha_1,\ldots,\alpha_n$ be the conjugates of $\alpha$, with  $\alpha_1=\alpha$.  Since $\alpha\neq0$, we have $a_0\neq0$. Moreover, $|a_0| =\left|\alpha_1\alpha_2\cdots\alpha_n\right|\leq 1$. Since $a_0$ is a nonzero integer, it follows that $|a_0|=1$. Hence $\prod_{i=1}^n|\alpha_i|=1$.
		Since $|\alpha_i|\leq1$, we must have
		$|\alpha_i|=1$ for each $i$. In particular, $|\alpha|=1$. Since
		$\alpha$ is real, $\alpha\in\{-1,1\}$. 
	\end{proof}
	
	\begin{lemma}\label{trace_rank}
		Let $M$ be a symmetric matrix with rational entries. If all eigenvalues of $M$ are algebraic integers, then $\trace(M^2)\geq\rank(M)$. Equality holds if and only if every nonzero eigenvalue of $M$ is $1$ or $-1$. Moreover, if $2\in\spec(M)$, then $\trace(M^2)\geq\rank(M)+3$.
	\end{lemma}
	\begin{proof}
		Let $M$ be an $n\times n$ matrix, and set $r=\rank(M)$. If $r=0$, then the result is immediate. Assume $r>0$. Let $\ld_1,\ldots,\ld_r$ be the non-zero eigenvalues of $M$, counted with multiplicity. Then the characteristic polynomial of $M$ is given by $\chi_M(x)=x^{n-r}(x-\ld_1)\cdots(x-\ld_r)$. Since $M$ has rational entries, $\chi_M(x)\in\Ql[x]$, and hence $\prod_{i=1}^{r}\ld_i\in\Ql$. Since each $\ld_i$ is an algebraic integer, $\prod_{i=1}^{r}\ld_i$ is an algebraic integer, and therefore $\prod_{i=1}^{r}\ld_i\in\Zl$. Since $\ld_i\neq0$ for $1\leq i\leq r$, we have $\left|\prod_{i=1}^{r}\ld_i\right|\geq1$. Thus by the arithmetic-geometric mean inequality,
		\[\frac{1}{r}\sum_{i=1}^r\ld_i^2\geq\left(\prod_{i=1}^{r}\ld_i^2\right)^{1/r}=\left|\prod_{i=1}^{r}\ld_i\right|^{2/r}\geq1.\]
		Consequently, 
		\[\trace(M^2)=\sum_{i=1}^{r}\ld_i^2\geq r=\rank(M).\]
		Equality holds if and only if $\ld_i^2=1$ for each $i$, that is, every non-zero eigenvalue of $M$ is $1$ or $-1$.
		
		Now suppose $2\in\spec(M)$, and label the eigenvalues so that $\ld_1=2$. If $r=1$, then $\trace(M^2)=4=\rank(M)+3$. Let $r\ge2$ and set $\mu=\prod_{i=2}^{r}\ld_i$. Then $\mu=\frac12\prod_{i=1}^{r}\ld_i\in\Ql$, giving that $\mu\in\Zl\setminus\{0\}$. Hence $|\mu|\ge1$. The arithmetic--geometric mean inequality
		applied to $\ld_2^2,\ldots,\ld_r^2$ gives $\sum_{i=2}^{r}\ld_i^2\geq (r-1)|\mu|^{2/(r-1)}\geq\ r-1$.
		Therefore
		\[
		\trace(M^2)=4+\sum_{i=2}^{r}\ld_i^2\ \ge\ 4+(r-1)=\rank(M)+3. \qedhere
		\]
	\end{proof}
	\section{Valency restrictions for perfect state transfer}
	In this section, we prove that a Cayley graph exhibiting perfect state transfer must have valency $1$, $2$, or $4$. Throughout this section, we write $G=\cay(\Gamma,S)$,  where $\Gamma$ is a finite group, not necessarily abelian,  and $\emptyset\neq S\subseteq\Gamma\setminus\{1\}$ satisfies $S=S^{-1}$.
	\subsection{Odd valencies and valency at least eight}
	We first show that perfect state transfer cannot occur in Cayley graphs with odd valency greater than one or valency at least eight.
	\begin{lemma}\label{odd_cayley}
		Let $\cay(\Gamma,S)$ be a Cayley graph over a finite group $\Gamma$ with $|S|>1$. If $|S|$ is odd, then $\cay(\Gamma,S)$ does not exhibit perfect state transfer.
	\end{lemma}
	\begin{proof}
		Let $d:=|S|$ be odd with $d>1$ and suppose $\cay(\Gamma,S)$ exhibits perfect state transfer. Let $\mu$ be an eigenvalue of the discriminant $P$ of $\cay(\Gamma,S)$. We show that $\mu\in\{-1,0,1\}$. By Lemma~\ref{reg}, $d\mu$ is an eigenvalue of the adjacency matrix $A$, and hence $d\mu$ is an algebraic	integer. By Lemma~\ref{algebraic_pst}, $2\mu$ is also an algebraic integer. Since $d$ is odd,
		\[\mu=d\mu-\frac{d-1}{2}(2\mu).\]
		Hence $\mu$ is an algebraic integer. Since each algebraic conjugate $\eta$ of $\mu$ is again an	eigenvalue of $P$, it follows that $\eta\in[-1,1]$. Hence by Lemma~\ref{alg_kro}, $\mu\in\{-1,0,1\}$.	For every $\mu\in\{-1,0,1\}$, we have
		\[
		T_\tau(\mu)=
		\begin{cases}
			\mu & \tau\text{ odd},\\
			1 & \tau\equiv0\pmod4,\\
			2\mu^2-1 & \tau\equiv2\pmod4.
		\end{cases}
		\]
		Since $\spec(P)\subseteq\{-1,0,1\}$ and $P$ is real
		symmetric, the spectral theorem gives
		\[
		T_\tau(P)=
		\begin{cases}
			P & \tau\text{ odd},\\
			I & \tau\equiv0\pmod4,\\
			2P^2-I & \tau\equiv2\pmod4.
		\end{cases}
		\]
		Since $d>1$, $P$ is not a permutation matrix, whereas $I$ fixes every vertex. Moreover,	every diagonal entry of $2P^2-I$ equals $2/d-1$.  Since $d$ is odd,  $2/d-1\neq0$. Thus, in case $2P^2-I$ is a permutation matrix, it must have fixed points. Hence by Lemma~\ref{pst_cayley}, $\cay(\Gamma,S)$ does not exhibit perfect state transfer.
	\end{proof}

	\begin{lemma}\label{geq8}
		Let $\cay(\Gamma,S)$ be a Cayley graph over a finite group $\Gamma$. If $|S|\geq 8$, then $\cay(\Gamma,S)$ does not exhibit perfect state transfer.
	\end{lemma}
	\begin{proof}
		Let $G=\cay(\Gamma,S)$, $n=|\Gamma|$, $d=|S|$, and let $A$ and $P$ denote the adjacency matrix and the discriminant matrix of $G$, respectively. Set $B:=2P=\frac{2}{d}A$ and $m_0:=\dim\ker B$.
		
		Suppose, to the contrary, that $G$ exhibits perfect state transfer at time $\tau$. By Lemma~\ref{pst_cayley}, $T_\tau(P)$ is a fixed-point-free permutation matrix and $T_\tau(P)^2=I$. Therefore $n$ is even, and $T_\tau(P)$ has eigenvalues $1$ and $-1$ with equal multiplicities $n/2$. Hence $\dim\ker (T_\tau(P)-I)=\dim\ker(T_\tau(P)+I)=n/2$. 
		
		We claim that $m_0\leq n/2$. If $m_0=0$, this is immediate. Suppose $m_0>0$, and let $0\neq y\in\ker B$. Since $\ker B=\ker P$, we have $Py=0$, and hence $T_\tau(P)y=T_\tau(0)y$. Since $T_\tau(P)^2=I$, we have $T_\tau(0)=1$ or $T_\tau(0)=-1$. If $T_\tau(0)=1$, then $y\in\ker(T_\tau(P)-I)$, and if  $T_\tau(0)=-1$, then $y\in\ker(T_\tau(P)+I)$. Therefore $\ker B\subseteq\ker(T_\tau(P)-I)$ or $\ker B\subseteq\ker(T_\tau(P)+I)$. In either case, $m_0\leq n/2$.
		
		Since $G$ is $d$-regular, $d$ is an eigenvalue of $A$, and hence $2$ is an eigenvalue of $B$.  Moreover, by Lemma~\ref{algebraic_pst}, every eigenvalue of $B$ is an algebraic integer. Therefore by Lemma~\ref{trace_rank}, together with $\trace(A^2)=nd$ and $m_0\leq n/2$, we obtain
		\[\frac{4n}{d}=\frac{4}{d^2}\trace(A^2)=\trace(B^2)\geq\rank(B)+3=(n-m_0)+3\geq\frac n2+3.\]
		This implies that $d<8$, contradicting the assumption $d\geq8$. This completes the proof.
	\end{proof}
	\subsection{Block decomposition and trace identities}\label{block_decompose}
	In this section, we discuss certain $2\times2$ block decompositions  and  trace identities of a symmetric matrix that will be used in the later sections. We write $\mathbf O$ for the zero matrix and $\mathbf 1$ for the all-ones vector, with sizes understood from the context. 
	
	Let $G=\cay(\Gamma,S)$ be a $d$-regular Cayley graph that exhibits perfect state transfer at time $\tau$. Let $A$ and $P$ be the adjacency and discriminant matrices of $G$, respectively. Since $G$ is $d$-regular, $P=\frac1d A$. Set $B:=2P$.	By Lemma~\ref{pst_cayley}, $T_\tau(P)$ is a fixed-point-free involutory permutation matrix. Let $\sigma$ be the corresponding permutation of $\Gamma$, defined by $T_\tau(P)\mathbf{e}_x=\mathbf{e}_{\sigma(x)}$ for each $x\in\Gamma$. Then $\sigma$ has no fixed points and $\sigma^2$ is the identity. It gives that $|\Gamma|$ is even.  Write $|\Gamma|=2m$.	Then $\Gamma$ can be partitioned into $m$ disjoint pairs $\{x_1,y_1\},\ldots,\{x_m,y_m\}$, where $\sigma(x_i)=y_i$ and $\sigma(y_i)=x_i$ for $1\leq i\leq m$. With respect to the ordering $x_1,\ldots,x_m,y_1,\ldots,y_m$ of $\Gamma$, we have
	\[
	T_\tau(P)=
	\begin{pmatrix}
		\mathbf O&I\\
		I&\mathbf O
	\end{pmatrix}.
	\] 
	Since $T_\tau(P)$ is a polynomial in $A$, it commutes with $A$. Comparing blocks in $AT_\tau(P)=T_\tau(P)A$, we obtain
	\[
	A=
	\begin{pmatrix}
		C&D\\
		D&C
	\end{pmatrix},
	\]
	where $C$ and $D$ are symmetric $0$-$1$ matrices.
	
	For $x\in V(G)$, let $N_G(x)=\{y\in V(G): y\text{ is adjacent to }x\}$. We now claim that $|N_G(x)\cap N_G(\sigma(x))|$ is independent	of $x$. Let $f$ be an automorphism of $G$. Since the permutation matrix of $f$ commutes with $A$, it also commutes with $T_\tau(P)$. Hence $f(\sigma(x))=\sigma(f(x))$ for each $x\in\Gamma$. Consequently,
	\[
	f\bigl(N_G(x)\cap N_G(\sigma(x))\bigr)
	=N_G(f(x))\cap N_G(\sigma(f(x))).
	\]
	Since $f$ is bijective and $G$ is vertex-transitive, the number $|N_G(x)\cap N_G(\sigma(x))|$ is independent of $x$. Let $r=|N_G(x)\cap N_G(\sigma(x))|$ for each $x$.

	For each $i,j$, the vertex $x_j$ is a common neighbor of $x_i$ and $y_i$ precisely when $C_{ij}=D_{ij}=1$. Further, $x_j$ is a common neighbor of $x_i$ and $y_i$ if and only if $y_j$ is a common neighbor of $x_i$ and $y_i$.  Hence
	\begin{equation}\label{crb1204}
		r=2\sum_{j=1}^{m}C_{ij}D_{ij}\quad\text{for}\quad1\leq i\leq m.
	\end{equation}
	In particular, $r$ is even and $0\leq r\leq d$. Moreover, since the graph is $d$-regular,
	\begin{equation}\label{crb1205}
		\sum_{j=1}^{m}(C_{ij}+D_{ij})=d\quad\text{for}\quad1\leq i\leq m.
	\end{equation}
	
	The eigenspaces of $T_\tau(P)$ corresponding to
	the eigenvalues $1$ and $-1$ are, respectively,
	\[
	E_+=\left\{\binom{z}{z}:z\in\mathbb R^m\right\}
	\quad\text{and}\quad
	E_-=\left\{\binom{z}{-z}:z\in\mathbb R^m\right\}.
	\]
   Since
	\[
	B=\frac{2}{d}
	\begin{pmatrix}
		C&D\\
		D&C
	\end{pmatrix},
	\]
	both $E_+$ and $E_-$ are $B$-invariant. Indeed,
	for every $z\in\mathbb R^m$,
	\[
	B\binom{z}{\pm z}
	=
	\frac{2}{d}
	\binom{(C\pm D)z}{\pm(C\pm D)z}\in E_\pm.
	\]
	Hence under the natural identifications
	\[
	\varphi_\pm:E_\pm\longrightarrow\mathbb R^m
	\text{ defined by }
	\varphi_\pm\binom{z}{\pm z}=z,
	\]
	the restrictions of $B$, considering as an operator, to $E_+$ and $E_-$ are
	represented, respectively, by
	\[
	B_+=\frac{2}{d}(C+D)
	\quad\text{and}\quad
	B_-=\frac{2}{d}(C-D).
	\]
	Consequently, the spectrum of $B$ is the union of the spectra of $B_+$ and $B_-$. Since $P=\frac12B$ and $T_\tau(P)$ is a polynomial in $P$, the restrictions of
	$T_\tau(P)$ to $E_+$ and $E_-$ are represented
	by $T_\tau(\frac12B_+)$ and $T_\tau(\frac12B_-)$,
	respectively. Moreover, $T_\tau(P)$ acts as $I$ on
	$E_+$ and as $-I$ on $E_-$, we obtain
	\begin{equation}\label{crb1229}
		T_\tau\left(\frac12B_+\right)=I\quad \text{and}\quad T_\tau\left(\frac12B_-\right)=-I.
	\end{equation}

	Noting $C$ and $D$ are $0$-$1$ matrices, we have
	\[
	(C_{ij}\pm D_{ij})^2
	=C_{ij}+D_{ij}\pm2C_{ij}D_{ij},
	\]
	and
	\[
	|C_{ij}-D_{ij}|
	=C_{ij}+D_{ij}-2C_{ij}D_{ij}.
	\]
	Since $B_+$ and $B_-$ are symmetric, it follows
	from \eqref{crb1204} and \eqref{crb1205} that
	\begin{align}
		\trace(B_\pm^2)&=\frac{4}{d^2}\sum_{i=1}^{m}\sum_{j=1}^{m}(C_{ij}\pm D_{ij})^2=\frac{4m(d\pm r)}{d^2},\label{trace_B_pm}\\
		\sum_{j=1}^{m}|(B_-)_{ij}|
		&=\frac{2}{d}\sum_{j=1}^{m} (C_{ij}+D_{ij}-2C_{ij}D_{ij})=\frac{2(d-r)}{d}\quad\text{for }1\leq i\leq m,~~\text{and}
		\label{mod_B_-}\\
		B_+\mathbf1
		&=\frac{2}{d}(C+D)\mathbf1
		=2\mathbf1.
		\label{ev_B_+}
	\end{align}
    
    Note that \eqref{crb1229}--\eqref{ev_B_+} hold under the assumption that the Cayley graph exhibit perfect state transfer at time $\tau$.
	
	\subsection{The valency-six case}
	We now show that Cayley graphs of valency six do not exhibit perfect state transfer, using the block decomposition
	and trace identities from Subsection~\ref{block_decompose}.
	\begin{lemma}\label{eqal6}
		Let $G:=\cay(\Gamma,S)$ be a Cayley graph over a finite group $\Gamma$. If $|S|=6$, then $G$ does not exhibit perfect state transfer.
	\end{lemma}
	\begin{proof}
		Suppose that $G$ exhibits perfect state transfer at time $\tau$ and apply the preceding calculation for $d=6$.  By \eqref{trace_B_pm},
		\begin{equation}
			\trace(B_\pm^2)=\frac{m(6\pm r)}{9}.
		\end{equation}
		By Lemma~\ref{algebraic_pst}, every eigenvalue of $B$ is an algebraic integer, and hence so is every eigenvalue of $B_\pm$. Since $r\geq0$, $\trace(B_-^2)=m(6-r)/9\leq2m/3<m$. Thus $B_-$ is singular. Indeed, if $B_-$ is non-singular, then Lemma~\ref{trace_rank} would give $\trace(B_-^2)\geq m$, a contradiction. Therefore $0\in\spec(B_-)$, and hence by \eqref{crb1229}, we have $T_\tau(0)=-1$.
		
		We now claim that $B_+$ is nonsingular. If $B_+$ is singular,  then there exists $x\neq 0$ such that $(B_+)x=0$. Taking  $y=(x,x)^t$, we have $By=0$ and $T_\tau(P)y=y$. However, $By=0$ gives $T_\tau(P)y=T_\tau(0)y=-y$, contradicting that $y\neq0$. Hence $B_+$ is nonsingular.
		
		By \eqref{ev_B_+}, $2$ is an eigenvalue of $B_+$. Therefore by Lemma~\ref{trace_rank}, $m(6+r)/9=\trace(B_+^2)\geq m+3$, giving that $r\geq4$. By \eqref{mod_B_-}, for  $1\leq i\leq m$, we have 
		\[\sum_{j=1}^{m}\left|(B_-)_{ij}\right|=\frac{6-r}{3}\leq\frac23.\]
		It follows that each eigenvalue $\theta$ of $B_-$ satisfies $|\theta|\leq 2/3$. Let $k=\rank(B_-)$. We claim that $k=0$. If $k>0$, then $B_-$ has exactly $k$ nonzero eigenvalues $\theta_1,\ldots,\theta_k$. Hence
		\[\trace(B_-^2)=\sum_{j=1}^{k}\theta_j^2\leq \frac{4k}{9}<k,\]
		contradicting Lemma~\ref{trace_rank}. Therefore $k=0$, and hence $B_-=\mathbf O$. Thus $C=D$. Then by \eqref{crb1205}, $\sum_{i=1}^{m} C_{ij}=3$ for $1\leq j\leq m$. Since $C$ is a symmetric matrix with entries in $\{0,1\}$, $\trace(C^2)=3m$. 
		
		Since $B_{+}$ is nonsingular, $C$ is nonsingular, and hence $\rank(C)=m$. Let $\ld$ be an eigenvalue of $C$. Since $C$ is an integer matrix, $\ld$ is an algebraic integer. Also, $B_{+}=\frac{2}{3}C$ implies that $2\ld/3$ is an algebraic integer. Hence $\ld/3$ is an algebraic integer. Thus every eigenvalue of $\frac{1}{3}C$ is an algebraic integer. Therefore by Lemma~\ref{trace_rank},
		\[m=\rank\left(\frac13C\right)\leq \trace\left(\frac19C^2\right)=\frac m3,\]
		contradicting that $m>0$. Thus $G$ does not exhibit perfect state transfer.
	\end{proof}

	We are now ready to prove our first main result, which shows that a Cayley graph exhibiting perfect state transfer must have valency $1$, $2$, or $4$.
	\begin{proof}[\textbf{Proof of Theorem~\ref{main1}}]
		The result follows immediately from Lemmas~\ref{odd_cayley}, \ref{geq8}, and \ref{eqal6}.
	\end{proof}

	\section{Abelian Cayley graphs of valency four}
	In this section, we classify connected Cayley graphs over finite abelian groups that exhibit perfect state transfer. Throughout this section, we write $G:=\cay(\Gamma,S)$, where $\Gamma$ is a finite abelian group with identity element $0$, and
	$S\subseteq\Gamma\setminus\{0\}$ satisfies $S=-S$ and $\langle S\rangle = \Gamma$. 
	
	Let $\Gamma$ be a finite abelian group. For $x\in\Gamma$, define the \emph{translation matrix} $\mathcal{M}_x$ by 
	\[\mathcal M_x\eu=\ee_{x+u}.\]
	Then $\mathcal M_x$ is a permutation matrix, and for all $x,y\in\Gamma$, 
	$\mathcal M_x\mathcal M_y=\mathcal M_y\mathcal M_x=\mathcal M_{x+y}$. 
	Moreover, $(\mathcal M_x)^{-1}=\mathcal M_{-x}=\mathcal M_x^t$. If $A$ is the adjacency matrix of the Cayley graph $\cay(\Gamma,S)$, then $A=\sum_{s\in S} \mathcal{M}_s$. Consequently, since $\cay(\Gamma,S)$ is $|S|$-regular, its discriminant matrix is
	\begin{equation}\label{crb1637}
		P=\frac{1}{|S|}\sum_{s\in S} \mathcal{M}_s.
	\end{equation}
	Note that $A\mathcal{M}_x=\mathcal{M}_xA$, and hence $f(P)\mathcal{M}_x=\mathcal{M}_xf(P)$ for any polynomial $f$.
	\begin{lemma}\label{ndef_pst}
		Let $G:=\cay(\Gamma,S)$ be an abelian Cayley graph. Then $G$ exhibits perfect state transfer from the vertex $u$ to another vertex $v$ at time $\tau$ if and only if $T_\tau(P)=\mathcal{M}_{v-u}$.
	\end{lemma}
	\begin{proof}
		By Lemma~\ref{pst_vertex}, $G$ exhibits perfect state transfer from  $u$ to $v$ at time $\tau$ if and only if $T_\tau(P)\eu=\ev$.  If $T_\tau(P)=\mathcal M_{v-u}$, then it is clear that $T_\tau(P)\eu=\ev$. Conversely, if $T_\tau(P)\eu=\ev$, then  for any vertex $x$, we have 
		\[
		T_\tau(P)\ex=T_\tau(P)\mathcal{M}_{x-u}\eu=\mathcal{M}_{x-u}T_\tau(P)\eu=\mathcal{M}_{x-u}\ev=\ee_{x+v-u}=\mathcal{M}_{v-u}\ex.
		\] 
	  Thus $T_\tau(P)=\mathcal{M}_{v-u}$.
	\end{proof}
	In particular, $\cay(\Gamma,S)$ exhibits perfect state transfer from $u$ to $v$ if and only if it exhibits perfect state transfer from $0$ to $h$, where $h:=v-u$.
	\begin{lemma}\label{crb1830}
		Let $G:=\cay(\Gamma,S)$ be an abelian Cayley graph with $|S|=4$. If $G$ exhibits perfect state transfer from $0$ to $h$ at time $\tau$, then $\ord(h)=2$, $\tau\equiv2\pmod4$, and $S+h=S$.
	\end{lemma}
	\begin{proof}
		Suppose that $G$ exhibits perfect state transfer from $0$ to $h$ at time $\tau$. By Lemma~\ref{ndef_pst}, $T_\tau(P)=\mathcal M_h$. Thus the permutation $\sigma$ corresponding $T_\tau(P)$ satisfies $\sigma(x)=x+h$ for $x\in\Gamma$. Since $\sigma$ is a fixed-point-free involution, we have $2h=0$ and $h\neq0$. Hence $\ord(h)=2$. We retain the notations of Subsection~\ref{block_decompose} 	with $d=4$. In particular, $|\Gamma|=2m$ and the paired vertices $x_i$ and $y_i$ satisfy $y_i=x_i+h$ for $1\leq i\leq m$. Further $B=\frac{1}{2}A$ and $B_\pm=\frac12(C\pm D)$, and that $r$ is even with $0\leq r\leq 4$. The trace formula becomes
		\begin{equation}\label{trace_4}
			\trace(B_\pm^2)=\frac{m(4\pm r)}{4}~\quad\text{and}\quad\trace(B^2)=2m.
		\end{equation}
		Moreover, by Lemma~\ref{algebraic_pst}, each eigenvalue  of $B$, and hence of $B_\pm$, is an algebraic integer.
		
		Since $2$ is an eigenvalue of $B$, Lemma~\ref{trace_rank} gives 
		\[2m=\trace(B^2)\ge\rank B+3.\]
		Hence $\rank B\leq2m-3$. Thus $B$ is  singular, and consequently
		$0\in\spec(P)$. Since $T_\tau(P)^2=I$, we have $T_\tau(0)=\pm1$, and it follows that $\tau$ is even. 
		
		We now claim that	
		\begin{equation}\label{excluded}
			\spec(B_-)\cap\{-1,1\}=\emptyset. 
		\end{equation} 
		Indeed, if $1$ or $-1$ were an eigenvalue of $B_-$, then $1/2$ or $-1/2$, respectively, would be an eigenvalue of $\frac12B_-$. Since $\tau$ is even, $T_\tau(1/2)=T_\tau(-1/2)$ and $T_\tau(1/2)\in\{1,-1/2\}$, contradicting \eqref{crb1229}.

		Suppose now that $\tau\equiv0\pmod4$. Then $T_\tau(0)=1$. Hence $B_-$ must be nonsingular; otherwise, $0$ would be an eigenvalue of $\frac12B_-$, contradicting \eqref{crb1229}. Therefore by \eqref{trace_4} and Lemma~\ref{trace_rank},
		\[m\leq\trace(B_-^2)=\frac{m(4-r)}{4}\leq m.\]
		Thus equality holds throughout. Hence by Lemma~\ref{trace_rank}, each eigenvalue of $B_-$ is either $1$ or $-1$, contradicting \eqref{excluded}. Therefore $\tau\equiv2\pmod4$.

		Now $\tau\equiv2\pmod4$ gives that $T_\tau(0)=-1$. This together with $T_\tau(\frac12B_+)=I$ implies that $B_+$ is non-singular. Since $2\in\spec(B_+)$, Lemma~\ref{trace_rank} gives
		\[
		\frac{m(4+r)}4=\tr(B_+^2)\ge m+3.
		\]
		This implies that $r>0$, and hence $r\ge2$.
		
		By \eqref{mod_B_-}, for $1\leq i\leq m$,
		\[
		\sum_{j=1}^{m}|(B_-)_{ij}|=\frac{4-r}{2}\le1.
		\]
		Hence each eigenvalue $\theta$ of $B_-$ satisfies $|\theta|\leq1$. Since $B_-$ has rational entries, every algebraic conjugate of $\theta$ is also an eigenvalue of $B_-$, and hence has absolute value at most $1$. Moreover, $\theta$ is a real algebraic integer. Thus by Lemma~\ref{alg_kro},  
		\[
		\spec(B_-)\subseteq\{-1,0,1\}.
		\]
	 This, along with \eqref{excluded}, gives that $\spec(B_-)=\{0\}$. Thus $B_-=\mathbf O$, and hence $C=D$. Consequently,
		\[
		A=
		\begin{pmatrix}
			C&C\\
			C&C
		\end{pmatrix}.
		\]
		Therefore the paired vertices have identical
		neighbourhoods, that is, for each $x\in\Gamma$,
		\[
		N_G(x)=N_G(\sigma(x))=N_G(x+h).
		\]
		Taking $x=0$, we obtain
		\[
		S=N_G(0)=N_G(h)=S+h.
		\]
		This completes the proof.		
	\end{proof}
	\begin{lemma}\label{pst_4}
		Let $\Gamma$ be a finite abelian group, and let $a,h\in\Gamma$ with $\ord(h)=2$. Let  $S=\{a,-a,a+h,-a+h\}$ with  $\langle S\rangle=\Gamma$ and $|S|=4$, and set $m=|\Gamma|/2$. Then $\cay(\Gamma,S)$ exhibits perfect state transfer from  $u$ to $v$ at time $\tau$ if and only if $v=u+h$, $m$ divides $\tau$, and $\tau\equiv2\pmod4$.
	\end{lemma}
	\begin{proof}
		Define the subgroup $H:=\{0,h\}$ of $\Gamma$. Since $\Gamma=\langle S\rangle$, we have $\Gamma=\langle a, h\rangle$. Hence $\Gamma/H$ is cyclic with $\Gamma/H=\langle a+H\rangle$. Moreover, $|\Gamma/H|=m$ and $a+H$ has order $m$.
		
		Define 
		\[F_1=\frac12(I+\mh),\quad F_2=\frac12(I-\mh),\quad\text{and}\quad L=\frac12(\mathcal{M}_{a}+\mathcal{M}_{-a}).\]
		Since $\ord(h)=2$, we find $\mh^2=I$. One can verify that,  for $i\in\{1,2\}$,
		\begin{equation}\label{crb1715}
			F_i^2=F_i,\quad LF_i=F_iL,\quad F_1F_2=0, \quad\text{and}\quad F_1+F_2=I.
		\end{equation}
		By \eqref{crb1637},
		\begin{align*}
			P&=\frac{1}{4}\left(\mathcal M_a+\mathcal{M}_{-a}+\mathcal{M}_{a+h}+\mathcal{M}_{-a+h}\right)\\
			&=\frac14\left(\mathcal M_a+\mathcal M_{-a}\right)\left(I+\mh\right)\\
			&=LF_1.
		\end{align*}
		Thus for $j\geq1$, $P^j=L^jF_1$. Let $f(x)=c_0+c_1x+\cdots+c_kx^k$ be a polynomial. Then
		\begin{align}
			f(P)&=c_0I+\sum_{j=1}^kc_jL^jF_1\nonumber\\
			&=c_0(F_1+F_2)+\sum_{j=1}^kc_jL^jF_1\nonumber\\
			&=f(L)F_1+f(0)F_2\label{crb1705}.
		\end{align}
		We now claim that for $j\geq0$,
		\begin{equation}\label{cheb_4_v}
			T_j(L)=\frac{1}{2}\left(\mathcal M_{ja}+\mathcal M_{-ja}\right).
		\end{equation}
		We prove \eqref{cheb_4_v} by induction on $j$. The cases $j=0$ and $j=1$ follow from $T_0(L)=I$ and $T_1(L)=L$. Suppose the identity holds for up to $j$, where $j\geq1$. Using the Chebyshev recurrence, we obtain
		\begin{align*}
			T_{j+1}(L)
			&=2LT_j(L)-T_{j-1}(L)\\
			&=\frac12(\mathcal M_a+\mathcal M_{-a})
			(\mathcal M_{ja}+\mathcal M_{-ja})
			-\frac12(\mathcal M_{(j-1)a}
			+\mathcal M_{-(j-1)a})\\
			&=\frac12(\mathcal M_{(j+1)a}
			+\mathcal M_{-(j+1)a}).
		\end{align*}
		Thus \eqref{cheb_4_v} holds for all $j\geq0$.

		Combining \eqref{crb1705} and \eqref{cheb_4_v}, we have
		\begin{equation}\label{crb1716}
			T_j(P)=\frac{1}{2}\left(\mathcal{M}_{ja}+\mathcal{M}_{-ja}\right)F_1+T_j(0)F_2.
		\end{equation}

		Suppose $\cay(\Gamma,S)$ exhibits perfect state transfer from $u$ to $v$ at time $\tau$. By translation, it is enough to consider perfect state transfer from $0$ to $w$, where $w=v-u$. By Lemma~\ref{pst_vertex},
		\begin{equation}\label{crb1746}
			T_\tau(P)\mathbf{e}_0=\mathbf{e}_w.
		\end{equation}
		Applying $F_2$ to both sides of \eqref{crb1746}, together with \eqref{crb1716} and \eqref{crb1715}, we have $T_\tau(0)F_2 \mathbf{e}_0=F_2\mathbf{e}_w$. This implies that 
		\begin{equation}\label{crb1725}
			T_\tau(0)(\mathbf e_0-\mathbf e_h)=\mathbf e_w-\mathbf e_{w+h}.
		\end{equation}
		Therefore $T_\tau(0)\neq0$ and $\{0,h\}=\{w,w+h\}$. Since $w\neq0$, we must have $w=h$, and hence $v=u+h$.
		
		Putting $w=h$ in \eqref{crb1725}, we have $T_\tau(0)(\mathbf e_0-\mathbf e_h)=(\mathbf e_h-\mathbf e_0)$. Therefore $T_\tau(0)=-1$, equivalently, $\tau\equiv2\pmod4$. 
		
		Since $w=h$, applying $F_1$ to \eqref{crb1746} and using \eqref{cheb_4_v} yields
		\begin{equation*}
			\frac{1}{4}\left(\ee_{\tau a}+\ee_{\tau a+h}+\ee_{-\tau a}+\ee_{-\tau a+h}\right)=\frac12\left(\ee_0+\ee_h\right).
		\end{equation*}
		This equality holds precisely when $\tau a\in H$, that is, $\tau a+H=H$. Since the order of $a+H$ is $m$, we have $m$ divides $\tau$.

		Conversely suppose $v=u+h$, $m$ divides $\tau$, and $\tau\equiv2\pmod4$. Then $\tau a\in H$. Therefore 
		\[\mathcal M_{\tau a}F_1=F_1=\mathcal M_{-\tau a}F_1.\]
		Also, $T_\tau(0)=-1$. Equation~\eqref{crb1716} gives $T_\tau(P)=F_1-F_2=\mh$. Thus by Lemma~\ref{ndef_pst}, $\cay(\Gamma,S)$ exhibits perfect state transfer from  $u$ to $v$ at time $\tau$.
	\end{proof}

	\begin{proof}[\textbf{Proof of Theorem~\ref{main2}}]
			Let $\Gamma$ be a finite abelian group, and let $S\subseteq\Gamma\setminus\{0\}$ satisfy
		$S=-S$, $\langle S\rangle=\Gamma$ and $|S|=4$. Suppose $\cay(\Gamma,S)$ exhibits perfect state transfer at time $\tau$. By translation, it occurs from $0$ to some $h(\neq0)$. Lemma~\ref{crb1830} gives that $\ord(h)=2$, $\tau\equiv2\pmod4$, and $S+h=S$. Set $H:=\{0,h\}$. The element $h\notin S$. Indeed, if $h\in S$ then $0=h+h\in S+h=S$, a contradiction. Consequently, $S$ is the union of two distinct nonzero $H$-cosets. Let $S=(a+H)\cup (b+H)$.  Since $S=-S$ and $-H=H$, we have $(a+H)\cup(b+H)=(-a+H)\cup(-b+H)$.
		In particular, $-a+H\subseteq S$. Since any two cosets of $H$ are either equal or disjoint, and $S$ is the union of the distinct cosets	 $a+H$ and $b+H$, we must have $-a+H=a+H$ or $-a+H=b+H$.

		Suppose  $-a+H=a+H$.  Then we also have	 $-b+H=b+H$. Thus $a+H$ and	$b+H$ are distinct nonzero elements of order $2$ in $\Gamma/H$. Moreover, since $\langle S\rangle=\Gamma$, we have $\Gamma/H=\langle a+H,b+H\rangle=\{H,a+H,b+H, a+b+H\}$. Since $H=\{0,h\}$, we have $\Gamma=\{0,h,a,a+h,b,b+h,a+b,a+b+h\}$ and $S=\{a,a+h,b,b+h\}$. Set $U:=S$ and $V:=\Gamma\setminus S=\{0,h,a+b,a+b+h\}$. Then every vertex of $U$ is adjacent to every
		vertex of $V$, with no edges within either set. Hence  $\cay(\Gamma,S)$ isomorphic to $K_{4,4}$, which does not exhibit perfect state transfer. Indeed,  for $K_{4,4}$, $\spec(P)=\{-1,0,1\}$. Then the spectral theorem gives
		\[
		T_\tau(P)=
		\begin{cases}
			P&\tau\text{ odd},\\
			I&\tau\equiv0\pmod4,\\
			2P^2-I&\tau\equiv2\pmod4.
		\end{cases}
		\]
		Every column of $P$ has four entries equal to $1/4$, $I$ fixes every vertex, and $2P^2-I$ has diagonal entries $-1/2$. Thus  by Lemma\ref{pst_cayley}, $K_{4,4}$ does not exhibit perfect state transfer.

		Therefore $-a+H=b+H$, giving that $S=(a+H)\cup (-a+H)=\{a,-a,a+h,-a+h\}$. Therefore by Lemma~\ref{pst_4}, $m$ divides $\tau$, where $m=|\Gamma|/2$. Since $\tau\equiv2\pmod4$, $4$ cannot divide $m$.
		
		Conversely, suppose there exist $a,h\in\Gamma$ with $h$ of order $2$ such that $4$ does not divide $|\Gamma|/2$ and  $S=\{a,\,-a,\,a+h,\,-a+h\}$. Let $m=|\Gamma|/2$. Take $\tau=2m$ or $\tau=m$ according as $m$ is odd or $m\equiv2\pmod4$. Therefore by Lemma~\ref{pst_4}, $\cay(\Gamma,S)$ exhibits perfect state transfer from $u$ to $u+h$ at time $\tau$.
		
		Let $\tau_{\operatorname{min}}$ be the minimum time at which $\cay(\Gamma,S)$ exhibits perfect state transfer from $u$ to $u+h$. Since $m\divides\tau_{\operatorname{min}}$, write $\tau_{\operatorname{min}}=km$. If $m$ is odd, then $km\equiv2\pmod4$ holds exactly when $k\equiv2\pmod4$, for which the least positive $k$ is $2$. If $m\equiv2\pmod4$, then $km\equiv2\pmod4$ holds exactly when $k$ is odd, for which the least positive $k$ is $1$. Thus $\tau_{\operatorname{min}}=\lcm(2,m)$.
	\end{proof}

	\section{Complete Classification of Cayley Representations}
	
	The results in the preceding sections, characterize connected abelian Cayley graphs that exhibit perfect state transfer. However, the same graph may have different Cayley representations, even over nonisomorphic groups. Thus a classification of the underlying graphs does not determine all possible groups and connection sets. In this section, we classify, up to
	equivalence, all pairs $(\Gamma,S)$ for which $\cay(\Gamma,S)$ exhibits perfect state transfer. This completes our classification of perfect state transfer on abelian Cayley graphs.
	
		\begin{table}[h!]
		\centering
		\renewcommand{\arraystretch}{1.5}
		\setlength{\tabcolsep}{20pt}
		\begin{tabular}{@{}cllcc@{}}
			\toprule
			Case & $\Gamma$ & $S$ & $h$ & $\tau_{\operatorname{min}}$\\
			\midrule
			\textup{(i)}
			& $\Zl_2$
			& $\{1\}$
			& $1$
			& $1$\\
			
			\textup{(ii)}
			& $\Zl_2\times\Zl_2$
			& $\{(1,0),(0,1)\}$
			& $(1,1)$
			& $2$\\
			
			\textup{(iii)}
			& $\Zl_{2n}$
			& $\{\pm1\}$
			& $n$
			& $n$\\
			
			\textup{(iv)}
			& $\Zl_{2m},~4\nmid m$
			& $\{\pm1,m\pm1\}$
			& $m$
			& $\lcm(m,2)$\\
			
			\textup{(v)}
			& $\Zl_m\times\Zl_2,~m\equiv2\pmod4$
			& $\{(\pm1,0),(\pm1,1)\}$
			& $(0,1)$
			& $m$\\
			\bottomrule
		\end{tabular}
		\caption{Classification of abelian Cayley graphs exhibiting perfect state transfer.}
		\label{tab:pst_classification}
	\end{table}
	
	\begin{proof}[\textbf{Proof of Theorem~\ref{complete_classification}}]
		Suppose that $G=\cay(\Gamma,S)$ exhibits perfect state transfer.
		By Theorem~\ref{main1}, we have $|S|\in\{1,2,4\}$.
		
		Suppose first that $|S|=1$, and write $S=\{a\}$. Since $S=-S$ and $a\neq0$, we have $\ord(a)=2$. Further,  $\langle S\rangle=\Gamma$ gives $\Gamma=\langle a\rangle\cong\Zl_2$. The isomorphism $\psi:\Gamma\to\Zl_2$ sending $a$ to $1$ gives the pair in \textup{(i)}.
		
		Next, suppose that $|S|=2$. Since $S=-S$, either $S$ consists of two distinct involutions or $S=\{a,-a\}$ with $a\neq-a$.
		In the first case, write $S=\{a,b\}$. Since $\Gamma$ is abelian and generated by the distinct nonzero involutions
		$a$ and $b$, we have $\Gamma=\{0,a,b,a+b\}\cong\Zl_2\times\Zl_2$. The isomorphism $\psi:\Gamma\to\Zl_2\times\Zl_2$ sending $a$ to $(1,0)$ and $b$ to $(0,1)$ gives the pair in \textup{(ii)}.
		
		In the second case, $\Gamma=\langle a\rangle$ is cyclic of order $q$, where $q=\ord(a)$. By Lemma~\ref{pst_cayley}, $q$ is even, that is, $q=2n$ for some integer $n$. The isomorphism $\psi:\Gamma\to\Zl_{2n}$ sending $a$ to $1$ therefore gives the pair in \textup{(iii)}.
		
		Finally, suppose that $|S|=4$. By Theorem~\ref{main2}, there exist $a,h\in\Gamma$ such that $\ord(h)=2$, $S=\{a,-a,a+h,-a+h\}$ and $4$ does not divide $m$, where $m:=|\Gamma|/2$.
		Set $H:=\langle h\rangle=\{0,h\}$. Since $\langle S\rangle=\Gamma$, we have $\Gamma=\langle a,h\rangle$, and hence $\Gamma/H=\langle a+H\rangle$. Thus $a+H$ has order $m$. Therefore $ma\in H$, and hence either $ma=h$ or $ma=0$.
		
		If $ma=h$, then $\Gamma=\langle a,h\rangle=\langle a\rangle$. Thus $a$ has order $2m$, and the isomorphism
		$\psi:\Gamma\to\Zl_{2m}$ sending $a$ to $1$ gives the pair in \textup{(iv)}.
		
		Suppose now that $ma=0$ and $m$ is even. The map $\psi:\Zl_m\times\Zl_2\to\Gamma$ by $\psi(j,k)=ja+kh$ 
		is clearly a homomorphism. Since $\Gamma=\langle a,h\rangle$, $\psi$ is surjective and both groups have order
		$2m$. Hence $\psi$ is an isomorphism, and $\psi^{-1}(S)=\{(1,0),(-1,0),(1,1),(-1,1)\}$. Since $4$ does not divide $m$ and $m$ is even, $m\equiv2\pmod4$ and $m\geq6$. This gives the pair in \textup{(v)}.
		
		If $ma=0$ and $m$ is odd, set $g:=a+h$. Since $ma=0$ and $2h=0$, we have $mg=ma+mh=h$. Thus $h\in\langle g\rangle$ and $a=g-h\in\langle g\rangle$. Hence $\Gamma=\langle g\rangle$ is cyclic of order $2m$. Moreover, $S=\{g,-g,g+h,-g+h\}$. The isomorphism $\psi:\Gamma\to\Zl_{2m}$ sending $g$ to $1$ gives the  pair in \textup{(iv)}.

		Conversely, in \textup{(i)}, we have $P=\mathcal M_1$,	and hence $T_1(P)=\mathcal M_1$. Thus perfect state
		transfer occurs at time $1$.
		
		For \textup{(ii)}, let $a=(1,0)$ and $b=(0,1)$. Then
		\[
		P=\frac12(\mathcal M_a+\mathcal M_b).
		\]
		Since $\mathcal M_a^2=\mathcal M_b^2=I$ and	$\mathcal M_a\mathcal M_b=\mathcal M_{a+b}$, we obtain $T_2(P)=2P^2-I=\mathcal M_{a+b}$. Hence perfect state transfer occurs from $u$ to	$u+a+b$ at time $2$ for every $u\in\Gamma$. Since $T_1(P)=P$ is not a permutation matrix,	the minimum transfer time is $2$.
		
		For \textup{(iii)}, we have
		\[
		P=\frac12(\mathcal M_1+\mathcal M_{-1}).
		\]
		An induction gives
		\[
		T_t(P)=\frac12(\mathcal M_t+\mathcal M_{-t})
		\quad\text{for}\quad t\geq0.
		\]
		Since $n=-n$ in $\Zl_{2n}$, we obtain $T_n(P)=\mathcal M_n$.  Hence perfect state transfer occurs from $u$ to $u+n$
		at time $n$ for every $u\in\Zl_{2n}$. Moreover, by Lemma~\ref{ndef_pst}, perfect state transfer occur from $u$ to $u+n$ at time $\tau$ if and only if $T_\tau(P)=\mathcal M_n$. Since $\mathcal M_x$ is a permutation matrix and $T_\tau(P)=\frac12(\mathcal M_\tau+\mathcal M_{-\tau})$, we find that $T_\tau(P)=\mathcal M_n$ if and only if $\tau\equiv-\tau\pmod{2n}$. Hence the minimum transfer time is $n$.

		In \textup{(iv)}, take $a=1$ and $h=m$. In \textup{(v)}, take $a=(1,0)$ and $h=(0,1)$. In both the cases, $S=\{a,-a,a+h,-a+h\}$, $h$ has order $2$, and $4$ does not divide $m$. Therefore by Theorem~\ref{main2}, the corresponding Cayley graphs exhibit perfect state transfer.

		An equivalence preserves $|\Gamma|$, $|S|$, and whether $\Gamma$ is cyclic. In \textup{(i)}, $|S|=1$. In \textup{(ii)} and \textup{(iii)}, $|S|=2$, and $\Gamma$ is non-cyclic in \textup{(ii)} but cyclic in \textup{(iii)}. In \textup{(iv)} and \textup{(v)}, $|S|=4$, and $\Gamma$ is cyclic in \textup{(iv)} but not in \textup{(v)}. Hence all these pairs are pairwise inequivalent.
	\end{proof}

	In Table~\ref{tab:pst_classification}, we summarize the classification of abelian Cayley graphs exhibiting perfect state transfer along with the transfer vertices and minimum transfer times. In each case, perfect state transfer occurs from $u$ to $u+h$ for every $u\in\Gamma$, at the minimum time $\tau_{\operatorname{min}}$.

	\subsection*{Statement of AI use}
	During the preparation of this work, the authors used ChatGPT to suggest proof ideas for few results and to generate examples. All statements and proofs were written and verified by the authors, who take full responsibility for the content of the manuscript.

	\subsection*{Acknowledgments}
	Koushik Bhakta acknowledges the financial support provided by the Prime Minister’s Research Fellowship (PMRF-ID: 1903298) scheme of the Government of India. 
	
	%%%%%%%%%%%%%%%%%%%%%%%%%%%%%%%%%%%%%%%%%%%%%%%%%%%%%%%%%%%%%%%%%%%%%%%%%%%%%%%%%%%%%%%%%%%%%%%%%%%%%%%%%%%%%%%%

\end{document}